\documentclass[11pt]{article}
\usepackage{GAFA2009}
\usepackage[paperwidth=193mm,paperheight=260mm,tmargin=25mm,bmargin=25mm,lmargin=25mm,rmargin=28mm]{geometry}
\usepackage{enumerate}
\usepackage[dvips]{graphicx}
\usepackage{color}              
\usepackage{epsfig}

\DeclareMathOperator{\thin}{Thin}
\DeclareMathOperator{\Mod}{Mod}

\DeclareMathOperator{\diam}{diam}

\DeclareMathOperator{\Proj}{Proj}
\DeclareMathOperator{\injrad}{injrad}
\begin{document}

\setfirstpage{0001}

\allowdisplaybreaks
     \abovedisplayskip=2pt plus 2pt minus 2pt
     \abovedisplayshortskip=-2pt plus 1pt minus 1pt
     \belowdisplayskip=2pt plus 2pt minus 2pt
     \belowdisplayshortskip=2pt plus 1pt minus 1pt

\newcommand{\TITLE}
{{
Divergence of geodesics
}}

\newcommand{\TItLE}
{{
Divergence of geodesics in Teichm\"uller
space and the mapping class group
}}

\newcommand{\AUTHOR}{{M.\ Duchin and K.\ Rafi}}

\newcommand{\AUtHOR}
{{Moon Duchin and Kasra Rafi}}

\TitlE
\AuthoR


\newcounter{introthm}[section]
\renewcommand{\theintrothm}{\Alph{introthm}}

\newenvironment{introthm}[1][{}]{\refstepcounter{introthm}
      \smallbreak\noindent\begingroup%
      \textbf{Theorem
      \theintrothm{\rm{#1}}.}
      \endgroup\nobreak\slshape\ignorespaces}{\smallbreak}

\newcounter{theorem}[section]
\renewcommand{\thetheorem}{\thesection.\arabic{theorem}}

\newenvironment{theorem}[1][{}]{\refstepcounter{theorem}
      \smallbreak\noindent\begingroup%
      \textbf{Theorem
      \thetheorem{\rm{#1}}.}
      \endgroup\nobreak\slshape\ignorespaces}{\smallbreak}

\newenvironment{lemma}[1][{}]{\refstepcounter{theorem}
      \smallbreak\noindent\begingroup%
      \textsc{Lemma
      \thetheorem{\rm{#1}}.}
      \endgroup\nobreak\slshape\ignorespaces}{\smallbreak}

\newenvironment{prop}[1][{}]{\refstepcounter{theorem}
      \smallbreak\noindent\begingroup%
      \textsc{Proposition
      \thetheorem{\rm{#1}}.}
      \endgroup\nobreak\slshape\ignorespaces}{\smallbreak}

\newenvironment{definition}[1][{}]{\refstepcounter{theorem}
      \smallbreak\noindent\begingroup%
      \textsc{Definition
      \thetheorem{\rm{#1}}.}
      \endgroup\nobreak\slshape\ignorespaces}{\smallbreak}




\newcommand{\thmref}[1]{Theorem~\ref{#1}}
\newcommand{\propref}[1]{Proposition~\ref{#1}}
\newcommand{\secref}[1]{\S\ref{#1}}
\newcommand{\lemref}[1]{Lemma~\ref{#1}}
\newcommand{\corref}[1]{Corollary~\ref{#1}}
\newcommand{\figref}[1]{Fig.~\ref{#1}}
\newcommand{\exref}[1]{Example~\ref{#1}}
\newcommand{\remref}[1]{Remark~\ref{#1}}



\newcommand{\co}{\colon\thinspace}
\newcommand{\emul}{\overset{.}{\asymp}}
\newcommand{\gmul}{\overset{.}{\succ}}
\newcommand{\lmul}{\overset{.}{\prec}}
\newcommand{\eadd}{\overset{+}{\asymp}}
\newcommand{\gadd}{\overset{+}{\succ}}
\newcommand{\ladd}{\overset{+}{\prec}}

\newcommand{\CC}{{\mathcal C}}
\newcommand{\SSS}{{\mathcal S}}
\newcommand{\cX}{{\mathcal X}}
\newcommand{\cY}{{\mathcal Y}}
\newcommand{\cZ}{{\mathcal Z}}
\newcommand{\T}{{\mathcal T}}
\newcommand{\M}{{\mathcal M}}
\newcommand{\F}{\mathcal{PMF}}
\newcommand{\MCG}{\Mod}
\newcommand{\G}{{\mathcal G}}
\newcommand{\CS}{\CC\SSS}
\newcommand{\RR}{t}
\newcommand{\Prod}{\mathop{\rm Prod}(\Gamma)}
\newcommand{\cc}{{\sf c}}
\newcommand{\thresh}{{\sf K}}

\renewcommand{\div}{\mathop{\rm div}}
\newcommand{\dist}{{\rm dist}}

\newcommand{\x}{{\sf X}}
\newcommand{\y}{{\sf Y}}
\newcommand{\z}{{\sf Z}}

\newcommand{\Teich}{Teich\-m\"uller~}

\renewcommand{\SS}{\mathcal{S}}
\newcommand{\hyp}{{\mathbb H}}
\renewcommand{\AA}{\mathbf{A}}
\newcommand\N{{\mathbb N}}
\newcommand\C{{\mathbb C}}
\newcommand\D{{\mathbf D}}
\newcommand\Z{{\mathbb Z}}
\newcommand\R{{\mathbb R}}
\newcommand\Q{{\mathbb Q}}
\newcommand\E{{\mathbb E}}

\newcommand{\ep}{\epsilon}
\newcommand{\p}{\partial}
\newcommand{\dps}{\displaystyle}
\newcommand{\sm}{\setminus}
\newcommand{\from}{\. \colon}
\newcommand{\st}{\bigm|}


\newcommand{\ab}{{\overline \alpha}}
\newcommand{\bb}{{\overline \beta}}
\newcommand{\gb}{{\overline g}}
\newcommand{\hb}{{\overline h}}


\Abstract{
We show that both \Teich space (with the \Teich metric) and the mapping
class group (with a word metric) have geodesic divergence that is 
intermediate 
between the linear rate of flat spaces and the exponential rate of 
hyperbolic 
spaces. For every two geodesic rays in \Teich space, we find that their 
divergence is at most quadratic. Furthermore, this estimate is shown to 
be sharp 
via examples of pairs of rays with exactly quadratic divergence.  The 
same 
statements are true for geodesic rays in the mapping class group. 
We explicitly describe efficient paths ``near infinity'' in both spaces.
}

\section{Introduction} \label{sec:intro}

The volume of a ball in \Teich space grows exponentially
fast as a
function of
its radius, as in the case of hyperbolic space. In this paper, we show 
that
despite this, the ``circumference'' of the ball grows only 
quadratically.
To be precise, for two geodesic
rays $\gamma_1$ and $\gamma_2$ in a proper geodesic space $X$ with a 
common basepoint $x\in X$, let their {\em divergence} be the infimal 
length of 
all paths connecting $\gamma_1(t)$ to $\gamma_2(t)$
which maintain a distance at least $t$ from the basepoint:
$$
\div(\gamma_1,\gamma_2,t)=\dist_{X\setminus B_t(x)}\bigl(
\gamma_1(t),\gamma_2(t)\bigr)\,,
$$
where $B_t(x)$ is the open ball of radius $t$ about $x$.  Note that 
$\div$
may be infinite in this generality.

\begin{figure}[ht]
\begin{center}
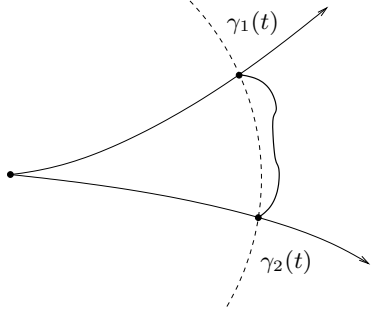
\begin{quote}
\caption{\small A path connecting $\gamma_1(t)$ to
$\gamma_2(t)$
that stays outside the open ball of radius $t$ centered at the 
basepoint.}
\end{quote}
\end{center}
\label{divergence}
\end{figure}

This definition is formulated to realize the trichotomy
that 
positive, zero, and negative curvature correspond to sublinear, linear, 
and superlinear divergence of geodesics, respectively; the 
correspondence is discussed further below.
The main results of this paper provide an upper 
bound 
for the divergence of geodesics in \Teich space and the mapping class 
group 
for surfaces of finite type.

Throughout this article, we compare non-negative functions $f(t)$ and 
$g(t)$
using the symbols
$\eadd,\emul,\asymp,\gadd,\gmul,\succ,\ladd,\lmul,\prec$
to denote equality or inequality with respect to an additive constant, a
multiplicative constant, or both, respectively, where the constants 
depend on
the topology of $S$ only. For example, $f(t) \asymp g(t)$ means that 
there are
constants $c_1$ and $c_2$ depending only on the topology of $S$ such 
that
$$
\tfrac 1{c_1}f(t)-c_2\leq g(t)\leq c_1f(t)+c_2\quad
\forall t\,.
$$
We may refer to this as the order of a function; for
instance if $f(t)\asymp t^2$ we may say that $f$ is on
the order of $t^2$.
Note that sometimes the distinction between $\ladd$ and $\lmul$ is 
crucial, as in many applications of Theorem
\ref{thm:product} below.
Note also that although much of the exposition is streamlined by this 
notation,
it is sometimes necessary to pay attention to the constants, as when 
considering functions $f(t)\asymp t^2$ and $g(t)\asymp t^2$ and trying 
to 
prove that $f(t)-g(t)\asymp t^2$.

With this equivalence relation, any two linear functions (respectively, 
polynomial 
of degree $n$) are identified.  Having $f\asymp 1$
means the function is bounded above.

\begin{introthm} \label{thm:upper}
Let $S$ be a surface of genus $g$ with $p$ punctures, such that 
$3g+p>4$.
Let $X$ be either the \Teich space $\T(S)$ with the \Teich metric or the
mapping class group $\MCG(S)$ with a word metric from a finite 
generating set.
For any pair of geodesic
rays $\gamma_1(t)$,$\gamma_2(t)$ with a common basepoint $x\in X$,
$$\div(\gamma_1,\gamma_2,t) \prec t^2.$$
\end{introthm}

To accomplish the quadratic upper bound, we explicitly construct paths 
that 
travel through chains of product regions. 
In the \Teich case, this path
travels through the {\em thin part} of $\T(S)$, which is stratified into 
regions 
which have a product structure, up to additive distortion.   Estimating 
the length of the path uses a combinatorial formula for
distance in \Teich space 
(\thmref{thm:distance}). 
In the case of the mapping class group, we use 
the 
quasi-isometrically embedded copies of $\Z^2$ generated by pairs of
Dehn twists about disjoint curves. 

\begin{figure}[ht]
\begin{center}
\includegraphics[width=2.2in]{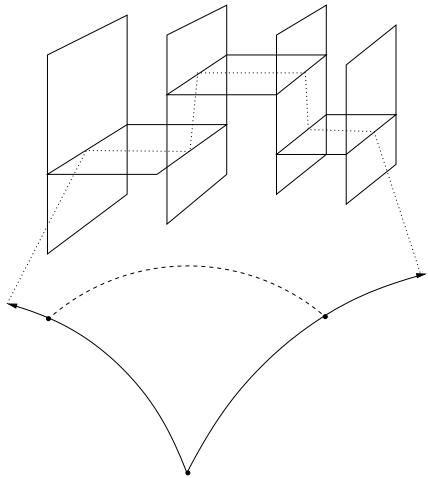}
\caption{\small An efficient path chains through product
regions in the space.}
\end{center}
\label{flatsfig}
\end{figure}

This theorem does not provide a 
quadratic or even linear lower bound for all divergence rates
in either space $X$, and indeed none is possible
since in both cases there are non-diverging pairs that pass every 
threshold 
separation. That is, based at every point in \Teich space or the mapping 
class group, and for arbitrarily large $M>0$,
there are pairs of rays with $\limsup \div(\gamma_1,\gamma_2,t)=M$.
One constructs these examples in \Teich space from pairs of quadratic
differentials with the same underlying topology
(see \cite{masur:CG}); in the mapping class group, one uses the many 
undistorted copies of $\Z^2$. 

On the other hand, the upper bound is sharp. 

\begin{introthm} \label{thm:lower}
For each $X$ as above, there exist pairs of rays
$\gamma_1,\gamma_2$ at every basepoint such that
$\div(\gamma_1,\gamma_2,t) \asymp t^2$.
\end{introthm}

Examples realizing the quadratic rate
are obtained from axes of pseudo-Anosov mapping classes. 
The lower bound on the divergence is furnished by means of 
{\em quasi-projection} theorems to these axes; in \Teich space, this is 
an
immediate application of a result of Minsky, while in the mapping class 
group
an analogous result is obtained using hyperbolicity of the curve 
complex.
(See ~\cite{behrstock:AG} for a different treatment of divergence in the 
case 
of the mapping class group.)

We remark that when $3g+p=4$ (that is, when $S$ is either
a 
once-punctured
torus or a four-times-punctured sphere), the \Teich metric on $\T(S)$
has negative curvature and the divergence is exponential. 
In
higher complexity, \Teich space is not hyperbolic (\cite{masur:NH}, 
\cite{minsky:PR}). 
However, there is a long-standing analogy
between the geometry of \Teich space and that of a hyperbolic space
(\cite{bers:EP}, \cite{kerckhoff:AG}, \cite{masur:TP},
\cite{wolpert:NC}, \cite{minsky:QP}). 
This paper provides another point 
of view
from which the \Teich space is quite different from a hyperbolic space.

\smallskip

These results should also be regarded as showing that $\T(S)$ and 
$\MCG(S)$
have {\em intermediate divergence} in the following sense.
A proper geodesic space $X$ is said to have {\em exponential divergence}
if there exists a threshold value  $D>0$ such that for all rays 
$\gamma_1,\gamma_2$,
$$\div(\gamma_1,\gamma_2,t_0)>D ~\hbox{for some}~ t_0 \implies
\div(\gamma_1,\gamma_2,t)\succ e^t.$$
Gromov hyperbolicity ($\delta$-hyperbolicity) is then equivalent to 
exponential 
divergence ($D=3\delta$ will suffice; see \cite{bridson:NPC}).
On the other hand, in flat spaces, 
every two geodesic rays have a linear divergence function.

A geodesic space can be said to have
{\em intermediate divergence} if no rays diverge faster than $f(t)$,
and some diverge at the rate $f(t)$, for a function growing 
super\-linearly
but subexponentially.  

In [G1,2], Gersten attributes to Gromov the 
expectation that there should be no nonpositively curved spaces of 
intermediate 
divergence. However, Gersten constructed such an example, giving a 
finite 
CAT(0) $2$-complex
whose universal cover possesses two geodesic rays which diverge
quadratically and such that no pair of rays diverges faster than
quadratically.  
Furthermore, in those papers, Gersten introduced a variation on 
divergence which
appeals directly to the metric with no reference to geodesics.
(He measures the lengths of paths staying far from a basepoint $x$
between pairs of points at a fixed distance from $x$,
with a slightly 
weaker 
equivalence relation on rates of growth than the one given here.)  
The 
advantage 
of this approach is that it produces a quasi-isometry invariant, so that 
the divergence 
of a finitely generated group can be discussed without specifying a 
generating set.  
Gersten and Kapovich--Leeb
(\cite{gersten:D3}, \cite{kapovich:3M})
study divergence in 3-manifold groups and find that some (namely 
graph-manifold 
groups) have quadratic divergence in the same sense as above.  
Here, we show that the same behavior occurs in the spaces $\T(S)$ and 
$\MCG(S)$,
which have been extensively studied in their own right.
Note that our Theorem A
implies quadratic divergence in Gersten's coarse sense, for the word 
metric with 
respect to any finite generating set.

\enlargethispage{-6pt}
By contrast, symmetric spaces have a gap in the possible
orders of
$\div(-,-,t)$ between linear and exponential rates; there, quadratic
divergence never occurs for any pair of rays  \cite{gromov:AI}.
Thus this paper also highlights the limitations of the long-standing 
analogies
between $\T(S)$ and symmetric spaces.

\subsection*{Acknowledgments.}
The authors would like to thank the Fields Institute for
hosting the 
workshop
at which these results were first conceived.  We also thank Howard 
Masur,
Yair Minsky, Saul Schleimer, and Dylan Thurston for very helpful
conversations during the conference.
In particular, Saul Schleimer observed that the
upper-bound arguments
for \Teich space would essentially carry over to the mapping class 
group.
Thanks also to Misha Kapovich, Rohit Thomas, and an anonymous referee
for helpful comments on the paper.


\section{Background}

\noindent
Let $S$ be an orientable, connected topological surface of genus $g$ 
with $p$ 
punctures. Throughout this paper, we assume that $3g+p>4$.
In this 
section, we 
review some background material and establish a few lemmas.

\subsection*{The space of curves and arcs.}
Let $Y$ be an essential subsurface of $S$ (if it is a proper subsurface, 
then it 
has boundary, and we write $\p Y$ for that set of curves).
By a {\it curve} in $Y$ we
mean a non-trivial, non-peripheral, simple closed curve in $Y$ and by an
{\it essential arc} we mean a simple arc, with endpoints
on the boundary \hbox{of $Y$}, that cannot be pushed to
the boundary of $Y$. 
In case $Y$ is not 
an annulus, the homotopy class of an arc is considered relative to the
boundary of $Y$; when $Y$ is an annulus, the homotopy class of an arc 
is considered relative to the endpoints of the arc.
We will use $D_\alpha$ to denote the Dehn twist about a simple closed 
curve $\alpha$.
\goodbreak

Let $\CC(Y)$ be the set of all homotopy classes of curves and
essential arcs on the surface $Y$.  
We define a distance on $\CC(Y)$ 
as follows: for $\alpha, \beta \in \CC(Y)$, define
$d_Y(\alpha,\beta)$ to be equal to one if $\alpha \neq \beta$
and $\alpha$ and $\beta$ can be represented by disjoint curves
or arcs. The metric on $\CC(Y)$ is the maximal metric having the
above property.  Thus, $d_Y(\alpha, \beta)= n$ if
$\alpha=\gamma_0, \gamma_1, \ldots, \gamma_n=\beta$ is some
shortest sequence of curves or arcs on $S$ such that successive
$\gamma_i$ are disjoint. Note that the notation $\CC(Y)$ 
is often used for the complex of curves on $Y$ (see \cite{minsky:CCI}
for definitions). 
Instead, we use $\CC(Y)$ to denote the zero-skeleton 
of the complex of curves and arcs, with distance induced
by the 
one-skeleton.  

In \cite{minsky:CCI}, Masur and Minsky show that $\CC(S)$
with the above metric is $\delta$-hyperbolic in the sense
of Gromov, and 
consequences
of this are discussed further below.
The curve complex is well known (and often proved) to be connected.
The isometry group of $\CC(S)$ is the mapping class group $\MCG(S)$ of 
$S$, 
as shown by Ivanov, Korkmaz, and Luo.

\subsection*{Subsurface projections.}
For $\alpha$ a curve in $S$, we define the {\it subsurface projection} 
of
$\alpha$ to the essential subsurface $Y$ as follows: Let
$$f:\bar S \to S$$
be a regular covering of $S$ such that $f_*(\pi_1(\bar S))$ is conjugate
to $\pi_1(Y)$ (this is called the $Y$-cover of $S$).
Since $S$ admits a
hyperbolic metric, $\bar S$ has a well-defined boundary at infinity, and 
we use the 
same notation $\bar S$ to denote the cover with its boundary added when 
appropriate. 
Let 
$\bar \alpha$ be the lift of $\alpha$ to $\bar S$. Components of $\bar 
\alpha$ 
that are essential arcs or curves on $\bar S$, if any, form a subset of 
$\CC(\bar S)$. The surface $\bar S$ is homeomorphic to $Y$. We call the 
corresponding subset of $\CC(Y)$ the {\em subsurface projection} of 
$\alpha$ to $Y$ 
and will denote it by $\alpha_Y$. If there are no essential arcs or 
curves in 
$\bar \alpha$, then $\alpha_Y$ is the empty set; otherwise we say that 
{\em $\alpha$ intersects $Y$ essentially}.
This projection only depends on $\alpha$ up to homotopy.

Let $\alpha$ and  $\alpha'$ be curves in $S$ that intersect a subsurface
$Y$ essentially. 
We define the {\em projection distance} between 
$\alpha$ and $\alpha'$ to be the maximum distance in
$\CC(Y)$ between the 
elements 
of the projections $\alpha_Y$ and $\alpha_Y'$, and denote it by 
$d_Y(\alpha,\alpha')$. If $Y$ is an annulus whose core is the curve 
$\gamma$, then we may denote $\CC(Y)$ by $\CC(\gamma)$ and 
$d_Y(\alpha,\alpha')$ by $d_\gamma(\alpha,\alpha')$.

Many of the results in the Masur--Minsky papers on the
geometry of the 
curve complex
proceed by analyzing the structure of $\CC(S)$ in terms of subsurfaces 
$Y$ by means 
of subsurface projection.  For instance, here is a useful result from 
\cite{minsky:CCI}.

\begin{theorem}[\ (Theorem 3.1 in \cite{minsky:CCI})]
\label{thm:sub-proj}
There is a constant $M_0$ such that if the projection
distance to some 
subsurface 
$Y$ satisfies 
$$
d_Y(\alpha, \beta)\geq M_0\,,
$$ 
then any geodesic in $\CC(S)$ connecting $\alpha$ to
$\beta$ intersects 
the 
$1$-neighborhood of $\partial Y$.
\end{theorem}

\subsection*{Markings.}
Following Thurston, a {\em pants decomposition} 
on $S$ is a 
maximal collection of disjoint curves;
these $3g-3+n$ curves are called the {\em pants curves} of that 
decomposition.
The {\em Fenchel--Nielsen coordinates} on $\T(S)$ are
obtained by 
assigning a length 
and a twist coordinate to
each curve in a pants decomposition.  Alternatively, a second set of 
curves may be 
chosen, transverse to the first,
so that lengths of all $6g-6+2n$ curves give approximate coordinates on 
$\T(S)$.
Along these lines, we define a {\it marking} on $S$ to be a set of pairs 
of curves
$\mu=\{ (\alpha_1, \beta_1), \ldots , (\alpha_m, \beta_m) \}$ such that
$\{ \alpha_1, \ldots , \alpha_m \}$ is a pants decomposition of $S$,
the curve $\beta_i$ is disjoint from $\alpha_j$ when  $i \not = j$,
and intersects $\alpha_i$ once (respectively, twice) if the surface
filled by $\alpha_i$
and $\beta_i$ is a once-punctured torus (respectively, a
four-times-punctured sphere).
The pants curves $\alpha_i$ are called the base curves of the marking
$\mu$.  For every
$i$, the corresponding $\beta_i$ is called the transverse curve to 
$\alpha_i$ in
$\mu$. When the distinction between the base curves and the transverse
curves is not important, we represent a marking as simply a set of $2m$
curves. Denote the space of all markings on $S$ by $\M(S)$.

Still following  \cite{minsky:CCII},
there are two types of {\em elementary moves} in $\M(S)$.
\begin{enumerate}
\item Twist: Replace $\beta_i$ by $\beta_i'$, where $\beta_i'$ is
obtained from $\beta_i$ by a Dehn twist or a half twist around
$\alpha_i$.
\item Flip: Replace the pair $(\alpha_i,\beta_i)$ with
$(\beta_i,\alpha_i)$;
also ``clean up'' by, for $j \not = i$, replacing $\beta_j$ with a curve
$\beta_j'$ that does not
intersect $\beta_i$ (the new base curve) in such a way that
$d_{\alpha_j}(\beta_j, \beta_j')$ is as small as possible
(see \cite{minsky:CCII} for details).
\end{enumerate}
In the first move, a twist can be positive or negative. A half twist
is allowed when $\alpha_i$ and $\beta_i$ intersect twice.
Masur--Minsky analyzed the geometry of the curve complex
by finding
efficient paths of markings through elementary moves.
The following theorem is a version of their result in which the high 
powers of Dehn 
twists are rearranged to appear consecutively; this adaptation
facilitates computations of \Teich distance through
changes of marking.

\begin{theorem}[\ \cite{rafi:CM}] \label{thm:path}
There exists a constant $\thresh$ (depending on $S$)
such that for any two markings $\eta_1,\eta_2$ on $S$,
there is a path of markings 
$$\eta_1=\mu_1,\ldots,\mu_n=\eta_2\,,
$$
where
$\mu_i$ and $\mu_{i+1}$ differ by an elementary move except that, 
for each $\alpha$ with $d_\alpha(\eta_1,\eta_2)\ge \thresh$, there is a 
unique index
$i_\alpha$ such that 
$$
\mu_{i_\alpha+1}=D_\alpha^p\mu_{i_\alpha}
\quad\text{with}\quad
|p|\eadd d_\alpha(\eta_1,\eta_2)\,.
$$
This path is efficient:
$$
n \asymp \sum_{Y\subseteq S} 
\big[d_Y(\eta_1,\eta_2)\big]_\thresh\,,
$$
where the sum is over non-annulus subsurfaces $Y$ and
$[N]_\thresh:=\begin{cases}N\,,&N\ge\thresh\,,\\0\,,&N<
\thresh\,.\end{cases}$
\label{markingchain}
\end{theorem}

\subsection*{Product regions in \Teich space.}
Let $\T(S)$ denote the \Teich space of $S$ equipped with the \Teich 
metric.
A point $\tau\in\T(S)$ is a hyperbolic metric on $S$ (constant curvature 
$-1$).
Minsky has shown that the thin part of \Teich space has a product-like
structure, as we now describe. Let $\Gamma$ be a set of disjoint curves
on $S$ and let $\thin_\ep(\Gamma)$ denote the set of points in \Teich 
space such that all curves from $\Gamma$ are short in
hyperbolic length:
$$
\thin_\ep(\Gamma):=\big\{\tau\in\T(S)\st l_\tau(\gamma)<
\ep ~\hbox{\rm for all}~ 
\gamma\in\Gamma \big\}\,.
$$
These cover the {\em $\ep$-thin part} of \Teich space,
which consists 
of
 those metrics with some short curve:
$$
\thin_\ep:=\big\{\tau\in\T(S)\st\injrad(\tau)<\ep\big\}\,.
$$
The {\em $\ep$-thick part} is the complement of
$\thin_\ep$.
Let $\Prod$ denote the product space,
$$
\Prod:=\T(S \setminus \Gamma) \times \prod_{\gamma \in \Gamma} 
\hyp_\gamma\,,
$$
where $S\setminus\Gamma$ is considered as a surface of
lower
complexity and each $\hyp_\gamma$ is a horoball in the copy of the 
hyperbolic plane parametrizing the Fenchel--Nielsen
coordinates corresponding to a short 
curve $\gamma$ (the $x$-coordinate in $\hyp_\gamma$
represents the twist parameter along $\gamma$ in $\sigma$
and the $y$-coordinate represents
the reciprocal of the length of $\gamma$ in $\sigma$; see 
\cite{minsky:PR}). Endow
$\Prod$ with the sup metric. Minsky has shown, for small
enough $\ep$, that $\thin_\ep(\Gamma)$ is well-approximated by
$\Prod$.

\begin{theorem}[\ (Product regions \cite{minsky:PR})]
\label{thm:product}
The Fenchel--Nielsen coordinates on $\T(S)$ give rise to
a natural homeomorphism $\pi \co \T(S) \to \Prod$.
There exists an $\ep_0> 0$ sufficiently small that this homeomorphism
restricted to $\thin_{\ep_0}(\Gamma)$ distorts distances
by a bounded
additive amount.
\end{theorem}

Note that $\T(S \setminus \Gamma) = \prod_Y \T(Y)$, where the product
is over all connected components $Y$ of $S \setminus \Gamma$.
Let $\pi_0$ denote the coordinate factor of $\pi$ mapping to
$\T(S \setminus \Gamma)$, let $\pi_Y$ denote the coordinate factor 
mapping
to $\T(Y)$, and, for $\gamma \in \Gamma$, let $\pi_\gamma$ denote the
coordinate factor mapping to $\hyp_\gamma$.

\subsection*{Short markings and \Teich distance.}
Let $\sigma$ be a point in the \Teich space $\T(S)$ of
$S$. A {\it short 
marking} 
\label{page:short} on $\sigma$ is a marking whose curves are chosen 
greedily 
to be as short as possible with respect to hyperbolic length.
That is, let $\alpha_1$ be a simple closed curve of minimal length in 
$\sigma$,
$\alpha_2$ a shortest curve disjoint from $\alpha_1$, and so on, to form
a pants decomposition of $S$ (the {\em Bers constant} gives an upper 
bound on the 
lengths of curves in a greedily chosen pants decomposition).  
Then, let the transverse curve $\beta_i$ be the
shortest curve intersecting $\alpha_i$ and disjoint from $\alpha_j$, $i 
\not = j$.
There are only finitely many choices in this process.
The following distance formula relates the \Teich
distance between two 
points
$\sigma_1$ and $\sigma_2$ to the combinatorics of short markings in
$\sigma_1$ and $\sigma_2$. Let $\ep_0$ be as before. Define
$\Gamma_{12}$ to be the set of curves that are
$\ep_0$-short
in both $\sigma_1$ and
$\sigma_2$, and, for $i=1,2$, define $\Gamma_i$ to be the
set of curves that are $\ep_0$-short in $\sigma_i$ but
not in 
$\sigma_{3-i}$.
Let $\eta_1$ and $\eta_2$ be short markings on $\sigma_1$
and $\sigma_2$, respectively.

\begin{theorem}[\ \cite{rafi:CM}] \label{thm:distance}
For sufficiently large $\thresh$,
the distance in $\T(S)$ between $\sigma_1$ and $\sigma_2$ is given
by the following formula:
\begin{multline} \label{eq:big}
  d_\T(\sigma_1, \sigma_2 ) \asymp
   \sum_Y \big[ d_Y(\eta_1, \eta_2) \big]_\thresh
+ \sum_{\alpha \not \in \Gamma_{12}} 
    \log \big[ d_\alpha(\eta_1,\eta_2)\big]_\thresh \\
+\max_{\alpha\in\Gamma_{12}}d_{\hyp_\alpha}(\sigma_1,
\sigma_2)
+
\max_{\begin{subarray}{l} \alpha \in \Gamma_i \\ i=1,2 \end{subarray}}
 \log \frac 1{l_{\sigma_i}(\alpha)}\,,
\end{multline}
where $\textstyle [N]_\thresh := \begin{cases}N\,,&N\ge 
\thresh\,,\\0\,,&N<\thresh\,,\end{cases}$ as before.
\end{theorem}

For the rest of the paper, we fix $\epsilon_0$ so that Minsky's
product regions theorem and \thmref{thm:distance} hold.

\subsection*{Subsurface distances and word length in the
mapping class group.} \label{subsec:sub-dist}
Consider the mapping class group $\MCG(S)$ of the surface $S$.
We fix a finite set $\AA$ of curves in $S$ such that the set
$\{D_\alpha :  \alpha\in\AA\}$ of Dehn twists around curves in $\AA$
generates $\MCG(S)$. (It follows that the curves in $\AA$ fill $S$.)
Equip $\MCG(S)$ with the word metric corresponding to this generating
set and denote the word length of an element $h \in \MCG(S)$ by $|h|$.

\begin{definition}[\ (Subsurface distance)]
Let $Y$ be a subsurface of $S$ and $\mu_1,\mu_2$ be two sets of curves
on $S$. We define
$$
d_Y(\mu_1,\mu_2) =
\max_{\alpha_1 \in \mu_1, \, \alpha_2 \in \mu_2} 
d_Y(\alpha_1, \alpha_2)\,,
$$
where $d_Y(\alpha_1, \alpha_2)$ is as defined above.
For $h \in \MCG(S)$, let
$$
d_Y(h) = d_Y \big( \AA, h(\AA) \big)\,,
$$
$\AA$ as above.
When $Y$ is an annulus whose core curve is $\gamma$, we may 
write $d_\gamma(h)$ instead of $d_Y(h)$.  
\end{definition}

Note that since $\AA$ fills $S$, $h(\AA)$ fills $S$ as
well
for every $h$; therefore, $h(\AA)$ intersects every subsurface
essentially and the projection of $h(\AA)$ to $Y$ is always
non-empty.

Next, we need the following theorem of Masur and Minsky which relates
the word length of a mapping class $h$ to the subsurface
distances
corresponding to $h$.  The sum in the statement of the theorem
is broken into two parts (summing over annular and non-annular 
subsurfaces) 
to highlight the comparison between 
this case and that in \thmref{thm:distance}.
The theorem essentially says that word length is comparable to the sum 
of the very large subsurface projections.

\begin{theorem}[\ (Word length in the mapping class group 
\cite{minsky:CCII})]
\label{thm:wordlength}
There is a constant  $\thresh_0$ such 
that for every threshold $\thresh\ge\thresh_0$ there exists $\cc$ 
such that for every $h\in\MCG(S)$,
\begin{equation}
\tfrac {|h|}\cc \leq \sum_\alpha\big[d_\alpha(h)\bigr]_
\thresh+\sum_Y\bigl[d_Y(h)\big]_\thresh  \leq 
\cc \, |h|\,,
\end{equation}
where the first sum is over all curves on $S$, the second sum is over
all
subsurfaces of $S$ that are not an annulus or a pair of pants and
$[N]_\thresh := \begin{cases}N\,,&N\ge \thresh\,,\\ 
0\,,&N<\thresh\,.\end{cases}$
\end{theorem}

We will sometimes refer to such a constant $\cc$ as the word-length 
constant
(for a threshold $\thresh$), taking
$\cc_0$ to be the value corresponding to $\thresh_0$.

\smallskip
We also need the following few simple lemmas.

\begin{lemma} \label{lem:triangle-inequality}
For $h_1, h_2,h_3 \in \MCG(S)$ and any subsurface
$Y$ of $S$, we have the following triangle inequality.
$$
d_Y\big(h_1(\AA),h_3(\AA)\big)\ladd d_Y\big(h_1(\AA),h_2(
\AA)\big)+d_Y\big(h_2(\AA),h_3(\AA)\big)\,.
$$
\end{lemma}

\pr 
The set $\AA$ has bounded diameter in $\CC(S)$, therefore its projection
to any subsurface $Y$ also has bounded diameter. The same is true for
$h_i(\AA)$, $i=1,2,3$. But $d_Y$ satisfies the triangle inequality in
$\CC(Y)$. Therefore, the above inequality is also true with the additive
error
of at most $3\diam(\AA)$ in $\CC(S)$.
\qed

One may think of $d_\alpha$ as measuring the relative twisting 
between two curves around $\alpha$.

\begin{lemma} \label{lem:twist}
For any curve $\alpha$ in $S$ and any mapping class $h\in\MCG(S)$,
$$d_\alpha(D_\alpha^n \, h) \gadd n- \cc \, |h|\,,
$$
where $\cc$ is the word-length constant.
\end{lemma}

\pr 
Using \lemref{lem:triangle-inequality}, we have
\begin{align*}
d_\alpha(D_\alpha^n \, h)  = d_\alpha \big(\AA, D_\alpha^n \, h(\AA) 
\big)
 & \gadd d_\alpha \big(h(\AA),  D_\alpha^n \, h(\AA) \big) 
   -  d_\alpha \big(\AA, h(\AA) \big) \\
 & \gadd n - d_\alpha(h)\,.
\end{align*}
Here, the last estimate is valid because $h(\AA)$ fills, so 
it contains a curve $\beta$ such that 
$d_\alpha(\beta,D_\alpha^n\beta)\eadd n$.
Finally, \thmref{thm:wordlength} implies that $d_\alpha(h) \ladd \cc \, 
|h|$
and the lemma follows.
\qed

\begin{lemma} \label{lem:invisible-twist}
If $\alpha$ and $\beta$ are disjoint, then 
$d_\alpha(D_\beta^n)$ has an upper bound that is independent of $n$.
\end{lemma}

\pr 
Let $\bar S$ be the annular cover of $S$ with respect to $\alpha$, 
let $\bar \alpha$ be the lift of $\alpha$ that is a closed curve and
$\bar \AA$ be the set of lifts of the elements of $\AA$ which 
intersect $\bar \alpha$.
Let $\gamma$ be a curve in $S$ intersecting $\alpha$ but not $\beta$ 
such 
that $d_\alpha(\AA, \gamma ) \asymp 1$ (the last condition can be 
obtained
after applying an appropriate power of $D_\alpha$ to $\gamma$). 
Let $\bar \gamma$ be a lift of $\gamma$ to $\bar S$ that intersects 
$\bar \alpha$.

Since $S$ is hyperbolic, $\bar S$ has a well-defined boundary at 
infinity. 
The lifts of $\beta$ are arcs with endpoints
on the boundary; a Dehn twist around $\beta$ in $S$ lifts to shearing 
along these
arcs in $\bar S$. 
Since $\bar \gamma$ is disjoint from all lifts of 
$\beta$,
no amount of shearing of $\bar \AA$ along lifts of $\beta$ can change 
the 
intersection number of $\bar \AA$ with $\bar \gamma$. That is,
$$
d_\alpha(\gamma,\AA)=d_\alpha(\gamma,D_\beta^n\AA)\asymp
1\,,
$$
and therefore
\begin{equation*}
d_\alpha(\AA, D_\beta^n \AA) \leq d_\alpha(\gamma,\AA)+d_
\alpha(\gamma,D_\beta^n\AA)\asymp
1\,.\tag*{$\scriptstyle\square$} 
\end{equation*}

\begin{lemma} \label{lem:high-twist}
For all words $h\in\MCG(S)$, all integers $n,m$, and all
curves $\alpha$ and $\beta$ such that $i(\alpha,\beta)=0$, we have
$$|D_\alpha^n D_\beta^m h| \gadd  \frac{\max(|n|,|m|)}\cc - |h|$$
for the word-length constant $\cc$.
\end{lemma}

\pr  Without loss of generality, $n\ge m\ge 0$. We have
\begin{align*}
|D_\alpha^nD_\beta^m \, h|
  &\gadd d_\alpha(D_\alpha^nD_\beta^mh)/\cc
       \tag{\thmref{thm:wordlength}}\\
  & = d_\alpha(D_\beta^nD_\alpha^m h)/\cc
       \tag{$D_\alpha$ and $D_\beta$ commute} \\
  & \eadd d_\alpha(D_\alpha^n \, h)/\cc
        \tag{Lemma~\ref{lem:invisible-twist}} \\
  & \gadd n/\cc -|h|\,.
        \tag{Lemma~\ref{lem:twist}} 
\end{align*}
This completes the proof.  \qed

Below, we will need to construct paths outside the ball
$B_t(e)$ in 
$\MCG(S)$.  
To initiate, it will be necessary to push away from the ball; the 
following lemma
says that for any Dehn twist, a  power of either it or its inverse 
accomplishes this.
Let 
$$\phi \colon \N \to \MCG(S)$$
be a geodesic ray in $\MCG(S)$ based at the identity.

\begin{lemma}[\ (Pushing off)] \label{lem:+or-}
There is a constant $d\in\N$ such that,
for any curve
$\alpha\in\AA$ and any nonnegative integers $n,m,\RR\in\N$, the 
magnitudes
$|D_\alpha^n \phi(d\RR)|$ and $|D_\alpha^{-m} \phi(d\RR)|$ are not both 
less than 
or equal to $\RR$.
\end{lemma}

\pr   
Let $c$ be the word-length constant from \thmref{thm:wordlength} and set 
$d>c+2$.
Let $h=\phi(d\RR)$ and suppose, for contradiction, that
$$|D_\alpha^n \,h|\le \RR \quad\text{and}\quad |D_\alpha^{-m}h|\le \RR$$
for some $m,n$.
Then
$|D_\alpha^{m+n}|\le 2\RR$.
But $d_\alpha(D_\alpha^{m+n}) \ge m+n-1$. 
Using the word-length constant from \thmref{thm:wordlength},
we obtain
\begin{equation} \label{eq:m+n}
m+n \leq 2\cc\, \RR +1\,.
\end{equation}
Since $D_\alpha$ is a generator for $\MCG(S)$,
$|D_\alpha^n| \leq n$.
Thus
$$
\RR\ge|D_\alpha^nh|\geq|h|-|D_\alpha^n|\geq d\RR-n\,.
$$
That is, $n  \ge (d-1)\RR$. Similarly, $m  \ge (d-1) \RR$ and
$$
m+n\ge 2(d-1)\RR>2(\cc+1)\RR=2\cc\, \RR+2\RR>2\cc\,\RR+1
\,,
$$
which is a contradiction.
\qed

\subsection*{Consequences of hyperbolicity.}
\label{subsec:hyp}

Recall that a space is called is {\em $\delta$-hyperbolic} for some 
$\delta>0$ if every geodesic triangle is $\delta$-thin: 
each side is 
contained in a $\delta$-neighborhood of the union of the
other two sides.  
Also, if a geodesic is regarded as an isometric embedding from (a 
subinterval of) $\R$ into $X$, then a {\em
$Q$-quasi-geodesic} replaces the equality 
with a coarse equality with additive and multiplicative
\hbox{constant $Q$}.
We collect here some standard consequences of hyperbolicity for use 
later 
in the paper. (See \cite{bridson:NPC} for a reference on hyperbolic 
spaces 
and quasi-geodesics.)

For the following three lemmas, let $X$ be a
$\delta$-hyperbolic space 
and let $L\subset X$ be a $Q$-quasi-geodesic line, ray,
or segment.
Suppose $a,b\in X$ and  $\bar a, \bar b \in L$ are such that $d(a,\bar 
a)$ and 
$d(b,\bar b)$ realize the distance from $a$ and $b$, respectively, to 
$L$.

\begin{lemma}[\ (Thin quadrilaterals)] \label{lem:quad}
There exist constants $M_1=M_1(\delta,Q)$ and $M_2=M_2(\delta,Q)$ such 
that 
if $d(\bar a,\bar b)>M_1$, then
any geodesic from $a$ to $b$
intersects the $M_2$-neighborhood of $L$.
\end{lemma}

\begin{lemma}[\ (Bounded shadows)]\label{lem:shadow}
There exists a constant $M_3=M_3(\delta,Q)$ such that 
for any geodesic segment $I$ from $a$ to $\bar a$, the closest-point 
projection 
of $I$ to $L$ has diameter $\le M_3$.
\end{lemma}

\begin{lemma}[\ (Bounded projection)]\label{lem:Lip}
There exists a constant $M_4=M_4(\delta,Q)$ such that
$d(\bar a,\bar b) \le M_4\cdot d(a,b) + M_4.$
\end{lemma}

Below, we will write $\Proj_L$ for the closest-point projection to a 
quasi-geodesic $L$.  
Since $\Proj_L(a)$ is of bounded diameter, this is a 
coarsely 
well-defined map.  

Recalling from above the Masur--Minsky result that
$\CC(S)$ is 
$\delta$-hyperbolic, 
we will reserve the notation $\delta$ for this particular hyperbolicity 
constant.


\section{Divergence in \Teich Space}
\noindent
As before, let $\T(S)$ denote the \Teich space of $S$ equipped with the 
\Teich metric. 
For a quadratic differential $q$ on a Riemann surface of topological 
type
$S$, let $[q]$ be the corresponding point of $\T(S)$. (That is, $[q]$ is 
the 
hyperbolic metric in the conformal class of $q$; for more background on 
\Teich space and quadratic differentials, see for instance 
\cite{strebel:QD},
\cite{abikoff:TT}, or \cite{hubbard:TT}.) Take $\sigma\in \T(S)$ and two 
quadratic differentials $q_1$ and $q_2$ such
that $[q_1]=[q_2]=\sigma$. Let $q_1(t)$ and $q_2(t)$ be the images of
$q_1$ and $q_2$, respectively, under the time-$t$ \Teich geodesic flow.
The maps
$$
t\mapsto\big [q_i(t)\big]\,, \quad i=1,2\,,
$$
from $[0,\infty)$ to $\T(S)$ are geodesic rays in $\T(S)$ emanating from
$\sigma$. We want to show that for all $\sigma,q_1,q_2$ as above and
$t>0$, there is a path from $[q_1(t)]$ to $[q_2(t)]$ in $\T(S)$
with length on the order of $t^2$ that stays outside
$B_t(\sigma)$, 
the open ball of radius $t$ in $\T(S)$ centered at $\sigma$ 
(see Fig.\ \ref{divergence}).

\enlargethispage{-3pt}
We will repeatedly use the same argument to show that a point 
in $\T(S)$ maintains a distance at least $t$ from
$\sigma$:  we show that 
the value $\epsilon(t)=e^{-2t}$ is small enough that any point in 
$\tau\in\T(S)$ with an $\epsilon$-short curve satisfies 
$d_\T(\tau,\sigma)\ge t$.
Then we carry out the appropriate sequence of elementary moves while 
maintaining some $\epsilon$-short curve at all times.

\subsection*{Constructing a path for the upper bound.}
To begin the progress from one ray to the other, we push off from the 
ball 
$B_t(\sigma)$ of radius $t$ around $\sigma$ so that the subsequent moves 
are guaranteed to stay far from $\sigma$. It suffices to construct, for 
sufficiently large $t$, a path between $\sigma_1=[q_1(3t)]$ and  
$\sigma_2=[q_2(3t)]$ whose length is of order $t^2$, while controlling 
the distance from $\sigma$.
The path will follow a sequence of elementary
moves, maintaining a sufficiently short curve at all times in order to 
ensure 
that we stay outside $B_t(\sigma)$.

Fix $\ep=e^{-2t}$ and suppose $\alpha$ on $S$ and $\tau\in\T(S)$ satisfy
$l_\tau(\alpha)\le \ep$.   In \cite{wolpert:LS}, it is shown that for a 
$K$-quasi-conformal map between Riemann surfaces, the
hyperbolic 
lengths of curves are changed by at most a multiplicative factor of $K$.  
Consequently the \Teich distance is bounded below by the
ratio of 
hyperbolic 
lengths for any particular curve.
$$
d_\T(\sigma, \tau) \geq \frac 12 \log 
\frac{l_\sigma(\alpha)}{l_\tau(\alpha)} \geq
\frac 12\log \frac {\ep_1}\ep\ge t\,.
$$
That is, this value of $\ep$ has the property described above
that if any curve on a surface $\tau$ is $\ep$-short, then $\tau\not\in 
B_t(\sigma)$.

We first give the argument under the assumption that the path starts and
ends in the $\ep_0$-thick part (that is, $\sigma_1,\sigma_2\not\in 
\thin_{\ep_0}$); 
we will treat the general case last. Recall that $\ep_0$ is chosen as in 
the 
product regions theorem.

Let $\mu_1, \ldots, \mu_n$ be the sequence
of markings described in \thmref{thm:path} such that $\eta_1=\mu_1$ is
a short marking on $\sigma_1$ and $\eta_2=\mu_n$ is a short
marking on $\sigma_2$. Note that the condition that $\sigma_1$ and
$\sigma_2$ are in the $\ep_0$-thick part of \Teich space implies
that the length of $\eta_i$ in $\sigma_i$ is bounded independent of $t$.

\goodbreak
For an elementary move on a marking, let the {\it
associated base
curve} be either the twisting curve, if the move is a twist, or the 
(initial) base curve in
the flipped pair, if the move is a flip. For the sequence of markings 
$\{\mu_i\}$,
let $\alpha_i$ be the associated base curve for the move $\mu_i\to 
\mu_{i+1}$
and let $\gamma_i$ be any base curve of the marking $\mu_i$ which is 
different
from $\alpha_i$ (the complexity condition $3g+p>4$ implies that there 
are at least
two base curves in every marking).  Note that $\gamma_i$ and 
$\gamma_{i+1}$
have intersection number zero. For any marking $\mu_i$ and $M$ 
sufficiently large, the set
$$
B_i=\big\{\tau\in\T(S):l_\tau(\mu_i)\le M\big\}\subset\T(
S)
$$
is nonempty and has bounded diameter 
(this follows, again, from the Wolpert formula). Fix such an $M$, which 
is suppressed in the notation $B_i$ from here on.

Since we can assume $t$ is sufficiently large, we will 
proceed taking 
$t{>}\max\{\diam B_i\}$ and $t>\log{1}/{\ep_1}$,
where $\ep_1$ is the injectivity radius of $\sigma$.

Let $\tau^1=\sigma_1$, $\tau^n= \sigma_2$ and, for
$1<i<n$, choose $\tau^i$ to be a point in $B_i$.

Let $\tau_{\gamma_i}^i$ be the point of $\T(S)$ that has the same
Fenchel--Nielsen coordinates as $\tau^i$ except that
$\gamma_i$
(which has bounded length in $\tau^i$) is pinched to have length
equal to $\ep$. 
Let $\tau^i_{\gamma_{i-1}}$ and $\tau^i_{\gamma_{i-1} \gamma_i}$
be defined similarly. Consider the piecewise geodesic path $P$ in 
$\T(S)$
defined by connecting up the points in the sequence below:
\begin{equation*}
\begin{split}
\sigma_1= \tau^1
\stackrel{p_1}\longrightarrow
\tau^1_{\gamma_1}
\stackrel{e_1}\longrightarrow
\tau^2_{\gamma_1}
\stackrel{p_2}\longrightarrow
\tau^2_{\gamma_1\gamma_2}
\stackrel{\bar{p}_1}\longrightarrow
\tau^2_{\gamma_2}
\stackrel{e_2}\longrightarrow
\tau^3_{\gamma_2}
\stackrel{p_3}\longrightarrow
\tau^3_{\gamma_2 \gamma_3}
\stackrel{\bar{p}_2} \longrightarrow
\tau^3_{\gamma_3}
\longrightarrow \cdots \qquad \\
\cdots \longrightarrow
\tau^{n-1}_{\gamma_{n-2} \gamma_{n-1}}
\stackrel{\bar{p}_{n-2}}\longrightarrow
\tau^{n-1}_{\gamma_{n-1}}
\stackrel{e_{n-1}}\longrightarrow
\tau^n_{\gamma_{n-1}}
\stackrel{\bar{p}_{n-1}}\longrightarrow
\tau^n = \sigma_2\,.
\end{split}
\end{equation*}
There are three kinds of steps:
\begin{enumerate}
\item[$(p_i)$]
\emph{pinches} the next curve, $\gamma_i$, by shortening
it to length $\ep$;
\item[$(e_i)$]
applies an \emph{elementary move} while keeping the curve
$\gamma_i$ short;
\item[$(\bar{p}_i)$]
\emph{releases} the previous curve, $\gamma_i$, by
restoring it to its
length before pinching.
\end{enumerate}

We may think of $p_i$, $\bar p_i$ and $e_i$ as paths in $\T(S)$ which 
are subpaths of $P$. 
We denote
the lengths of these paths by $|p_i|$, $|\bar p_i|$ and $|e_i|$ 
respectively.
We will bound the lengths of these subpaths and show that
they stay 
outside $B_t(\sigma)$
in order to complete the proof of Theorem A for $\T(S)$.

Except along $p_1$ and $\bar{p}_{n-1}$, at every point in the path $P$, 
at least one curve has length equal to $\ep$; therefore, by the choice 
of
$\ep$, these points are outside  $B_t(\sigma)$.
The lengths of the subpaths $p_1$ and 
$\bar{p}_{n-1}$ are (up to an additive error) equal to $t$ and they have 
one point at distance $3t$ from $\sigma$; therefore these
subpaths are 
also \hbox{outside $B_t(\sigma)$}.

We will estimate the lengths of the $p_i$ and $\bar{p}_i$
by showing 
that
$$
d_\T(\tau^i, \tau^i_{\gamma_i}) \eadd
d_\T(\tau^i, \tau^i_{\gamma_{i-1}}) \emul t\,,
$$
and
$$
d_\T(\tau^i_{\gamma_{i-1}}, \tau^i_{ \gamma_{i-1} \gamma_i}) \eadd
d_\T(\tau^i_{\gamma_{i-1}\gamma_i},\tau^i_{\gamma_i})
\emul t\,,
$$
To estimate the distance from $\tau^i$ to
$\tau^i_{\gamma_i}$,
consider first pinching the curve $\gamma_i$ in $\tau^i$ to obtain a new 
metric
$\tau \in \T(S)$, where $\gamma_i$ has  length  $\ep_0$.
The path from $\tau^i$ to $\tau$ has bounded length because 
both $\tau^i$ and $\tau$ are in $B_i$, which has bounded diameter. 
Now letting $\Gamma = \{ \gamma_i\}$, we note that $\tau, 
\tau^i_{\gamma_i}$
are both in $\thin_{\ep_0}(\Gamma)$ and \thmref{thm:product} applies.
But their projections to $S \setminus \Gamma$ are identical. Therefore, 
their distance, up to an additive error, is equal to the
distance in 
$\hyp_{\gamma_i}$
between $\pi_{\gamma_i}(\tau^i_{\gamma_i})$ and $\pi_{\gamma_i}(\tau)$, 
their projections to that factor.
The points $\tau$ and $\tau^i_{\gamma_i}$ 
have the same twisting parameters around $\gamma_i$ (the 
$x$-coordinates 
of their projections to $\hyp_{\gamma_i}$ are the same), but different 
length 
parameters (the $y$-coordinates are $1/\ep_0$ and $1/\ep$ 
respectively). 
Therefore
$$
d_\T(\tau^i_{\gamma_i}, \tau) \eadd
d_{\hyp_{\gamma_i}}\big(\pi_{\gamma_i}(\tau^i_{\gamma_i}), 
\pi_{\gamma_i}(\tau)\big)
= \tfrac 12 \log \tfrac{\ep_0}{\ep} \eadd t\,.
$$

Next, we estimate the lengths of the $e_i$ by showing that
$$
d_\T(\tau^i_{\gamma_i}, \tau^{i+1}_{\gamma_i}) \eadd
\begin{cases}
1&\text{if $\mu_i$ and $\mu_{i+1}$ differ by an elementary
move\,,}\\
\log p&\text{if $\mu_{i_\alpha+1}=D_\alpha^p\mu_{i_\alpha
}$\,, \quad 
$|p|\eadd
d_\alpha(\eta_1,\eta_2)\ge \thresh$\,.}
\end{cases}
$$
 If the elementary move from
$\mu_i$ to $\mu_{i+1}$ is a simple twist or flip, then the length of the
segment $(e_i)$ is bounded. This is because 
$\pi_{\gamma_i}(\tau^i_{\gamma_i})$ and 
$\pi_{\gamma_i}(\tau^{i+1}_{\gamma_i})$
are within bounded distance, and therefore
\begin{align*}
d_\T( \tau^i_{\gamma_i}, \tau^{i+1}_{\gamma_i}) 
& \eadd d_{S \setminus \gamma_i}
  \left( \pi_0( \tau^i_{\gamma_i}) , \pi_0(\tau^{i+1}_{\gamma_i}) 
\right)\\ 
& \eadd d_{S \setminus \gamma_i} \left( \pi_0( \tau^i) , 
\pi_0(\tau^{i+1}) \right)\\
&\ladd d_\T( \tau^i, \tau^{i+1}) \asymp 1\,.
\end{align*}
But there can also be high powers of twists:
for every $\alpha$ where $d_\alpha(\mu_1, \mu_n)$ is large
there is an index $i_\alpha$ where
$\mu_{i_\alpha+1} = D_\alpha^p \mu_{i_\alpha}$, as in \thmref{thm:path}.
Here the length of $(e_i)$ is on the order of
$\log d_\alpha(\mu_1, \mu_n)$. Therefore the total length of
$P$ is
\begin{equation*}
 \sum_i |p_i| + |\bar p_i| + |e_i| \asymp n\,t+\sum_{
\alpha}\log\big[d_\alpha(\mu_1,\mu_n)\big]_\thresh\,.
\end{equation*}
But \thmref{thm:path} and \thmref{thm:distance} imply
$$
n\asymp\sum_{Y\subseteq S}\big[d_Y(\mu_1,\mu_n)\big]_
\thresh \leq
d_\T(\sigma_1, \sigma_2) \leq 6t \,.
$$
It follows that the length $|P|$ is at most on the order
of 
$$ (6t) \, t + 6t \emul t^2.$$
This tells us in particular that the elementary moves (which move within
a single product region) contribute
negligibly to the total length of the path when the quadratic rate is 
realized; 
the pinch-and-release steps 
(which pivot from one product region to the next) account
for the whole length of the path, asymptotically.

\smallskip
Above, we assumed that the starting and ending points of
the path were in the $\ep_0$-thick part.
This condition was used in the definition
of $\tau^i$ (for example we have assumed that $\tau^1 \in B_1$).
This assumption is not always true. However, we can modify the beginning
and the end of the path above to accommodate the general case as 
follows:

Let $\tau$ be a point in $\T(S)$ with the same
Fenchel--Nielsen 
coordinates
as $\sigma_1$ except that the length of $\gamma_1$ is $\ep$.  Let
$\tau^1_{\gamma_1}$ be the point obtained by increasing the lengths of 
other
short curves (lengths less than $\ep_0$) to a moderate length.
We will show that these paths 
$P_1=[\sigma_1,\tau]$ and $P_2=[\tau,\tau^1_{\gamma_1}]$
are outside $B_t(\sigma)$ and have length of \hbox{order
$t$};
the rest of the 
calculation
reverts to the arguments above.

Note that there is a lower bound on the
$\sigma_1$-lengths of all 
curves.
That is, again using \cite{wolpert:LS} and recalling that 
$\ep_1$ is the injectivity radius at $\sigma$,
$$
\frac 12 \log \frac{ l_\sigma(\gamma)}{l_{\sigma_1}(\gamma)} \leq 3t
\quad \forall \gamma
\quad\Longrightarrow\quad
l_{\sigma_1}(\gamma) \ge \frac{l_\sigma(\gamma)}{e^{6t}}\ge \ep_1 
e^{-6t}. 
$$
Since there is no twisting or other change in the marking, it
follows from \thmref{thm:distance} that  both of the
paths $P_1$ and $P_2$ have length at most on the order of
$t$ (the only contribution to 
distance comes
from the length ratios).
The path $P_2$ stays outside $B_t(\sigma)$
because along this path the curve $\gamma_1$ has length $\ep$. To see 
that
the path $P_1$ is outside $B_t(\sigma)$ we take two cases.
If the length of $\gamma_1$ in $\sigma_1$ is less than $\ep$, then it 
remains less
than $\ep$ along this path and we are done. If
$\ep_0 \geq l_{\sigma_1}(\gamma_1) \geq \ep$, then the length of
this path is no more than $\frac 12 \log \frac{\ep_0}{\ep} \eadd t$
by the product regions theorem (\thmref{thm:product}).
 But $d_\T(\sigma,\sigma_1)=3t$, which completes the argument.

Applying the same modifications to the end of path $P$, we obtain the
desired result for the general case.

\subsection*{Example realizing the quadratic rate.}
To see that the quadratic estimate is sharp, we must furnish an example
of a pair of rays whose divergence rate is exactly quadratic.  We will 
use
a {\em quasi-projection} result of Minsky \cite{minsky:QP}
which shows 
that intervals which are
far (relative to their length) from a cobounded geodesic (segment, ray, 
or line)
project to sets whose diameter is uniformly bounded above. Recall that 
an
$\epsilon$-cobounded geodesic in $\T(S)$ is one which stays
in the $\epsilon$-thick part.  
(Note that part of the proof entails 
that $\Proj$, the closest-point projection, is coarsely
well-defined in the 
following setting.)

\begin{theorem}[\ (Quasi-projection for $\T(S)$ 
\cite{minsky:QP})] 
\label{thm:QP}
For every $\ep>0$ there are constants $b_1,b_2$ depending on $\epsilon$ 
and
the topology of $S$ such that the following holds. 
Let $G$ be an 
$\epsilon$-cobounded geodesic in $\T(S)$, suppose $\tau\in\T(S)$ 
satisfies $d(\tau,G)>b_1$,
and let  $r=d(\tau,G)-b_1$.  Then  
$$
\diam\bigl(\Proj_G (B_r(\tau))\big)\le b_2\,.
$$
\end{theorem}

It is straightforward to replace
the geodesic in the statement of the theorem with a quasi-geodesic.
Also note that it suffices to check a bounded geometry condition for the 
endpoints of the geodesic (in $\T(S)$ or on the Thurston boundary) to 
ensure that it stays in the thick part 
~\cite{rafi:thesis}.

To apply this theorem, we may for instance choose $q$ and $q'=-q$ to 
point in
opposite directions along the axis of a pseudo-Anosov mapping
class.  (Coboundedness is guaranteed because pseudo-Anosov axes project
to closed curves in moduli space.)

\begin{prop}
For any $\epsilon$-cobounded geodesic $G$ in $\T(S)$
and any point $\sigma \in G$, let $\sigma_t$ and $\sigma_t'$ be the
points on $G$ at distance $t$ from $\sigma$.  
Then for any path $P$ in $\T(S)$ 
from $\sigma_t$ to $\sigma_t'$ which maintains a distance
at least $t$ from 
$\sigma$, 
$$|P| \succ t^2.$$
\end{prop}

\pr 
Take $G_0$ to be the subsegment of $G$ of
length $t$, centered at $\sigma$.  Then any path $P$ between $\sigma_t$
and $\sigma_t'$ outside $B_t(\sigma)$ maintains a
distance at least $t/2$ 
from $G_0$ at all
times.

Now if $P$ is a path connecting $\sigma_t$ and $\sigma_t'$ in $\T(S)$ of
length $|P|$, we can divide it into pieces of length $t/2$, taking
$P=I_1\cup\cdots\cup I_n$ for successive pieces of length $t/2$
(with $I_n$ possibly shorter). The number of these intervals, $n$, 
satisfies
$$
2|P|\le nt \le 2|P| + t\,.
$$

Each $I_i$ maintains a distance at least $t/2$ from $G_0$.
Since $|I_i|\le t/2$ and $t$ is large, $I_i$ is covered
by two balls of radius $t/2-b_1$, centered on the endpoints of the 
interval.
Thus Minsky's quasi-projection theorem assures that 
$|\!\Proj_{G_0}(I_i)|<2b_2$ for each $i$.  
If the endpoints of $I_i$ are $x_{i-1}$ and $x_i$ and $\bar 
x_{i-1}$ and $\bar x_i$ are closest-point projections
\hbox{to $G$}, we see that $d(\bar 
x_{i-1},\bar x_{i})\le 2b_2$.
But then $t=d(\bar x_0, \bar x_n) \le 2nb_2$, and
since
$nt\asymp |P|$, we obtain $|P|\succ t^2$.
\qed


\section{The Mapping Class Group}

\subsection*{Constructing a path for the upper bound.}
Consider two infinite distinct geodesic rays, 
$$
\phi:\N\to\MCG(S)\quad\text{and}\quad \psi:\N\to\MCG(S)\,,
$$
in $\MCG(S)$ emanating from a common point (without loss of generality, 
the identity element) in $\MCG(S)$. We want to show that, for all $\RR 
\in \N$,
there is a path from $\phi(\RR)$ to $\psi(\RR)$ in $\MCG(S)$, such
that no point in this path is within distance $\RR$ of
the origin, and 
whose length is on the order of $\RR^2$.
We construct this path by traveling 
iteratively through 
chained copies of $\Z^2$, each copy generated by Dehn twists about a 
pair of 
disjoint curves.  However, first we need to move far enough from the 
identity, 
to points $\phi(d\RR)$ and $\psi(d\RR)$, so that the first few steps of 
the
sequence are
sure not to backtrack near the identity.  This is accomplished
by taking $d$ from \lemref{lem:+or-}.

Let $B_\RR=B_\RR(e)$ be the ball of radius $\RR$ in $\MCG(S)$
about the identity. The segments $[\phi(\RR), \phi(d\RR)]$ and
$[\psi(\RR), \psi(d\RR)]$ stay outside $B_\RR$ and their lengths are of
order $\RR$. To prove \thmref{thm:upper} for $\MCG(S)$, it is sufficient 
to build a 
path $P$ between the two points $u=\phi(d\RR)$ and $v=\psi(d\RR)$ that
stays outside $B_\RR$ and has length on the order \hbox{of
$\RR^2$}.

The path from $u$ to $v$ will involve high powers of Dehn
twists 
arranged
in ``switch moves" from one twist flat to the next.  We fix the exponent 
  $m=m(\RR)$ to be larger than
$(3d\cc +\cc)\RR$ (but of order $\RR$).

\begin{lemma}[\ (Switch moves)]\label{lem:gadget}
For any two curves $\alpha,\beta\in\AA$, any 
$w\in\MCG(S)$ with $|w|\le 3d\RR$, and $m$ chosen as above,
there is a path $(\star)$ from $D_\alpha^m w$ to
$D_\beta^m w$ staying outside $B_\RR$ and
of length $\prec \RR$.
\end{lemma}

\pr 
For the curves $\alpha,\beta$, we fix a chain of curves
\begin{equation}
\label{chain}
\{\gamma_i^{\alpha\beta}\}=\{\gamma_1^{\alpha\beta},
\ldots,\gamma_k^{\alpha\beta}\},
\quad\hbox{with}~ k=k(\alpha,\beta)\,,
\end{equation}
having the property that
$$\alpha-\gamma_1^{\alpha\beta}-\cdots-\gamma_k^{\alpha\beta}-\beta$$
is a path in the curve complex, that is, each adjacent pair of curves is
disjoint
on $S$. When there is no ambiguity about $\alpha$ and $\beta$, we denote
these curves simply by $\{\gamma_i\}$. Define $(\star)$ to be the 
following path 
from $D_\alpha^m w$ to $D_\beta^m w$:
\begin{equation}
D_\alpha^m \, w
\stackrel{D_{\gamma_1}^m}\longrightarrow
D_{\gamma_1}^m \, D_{\alpha}^m \, w
\stackrel{D_{\alpha}^{-m}}\longrightarrow
D_{\gamma_1}^m \, w
\stackrel{D_{\gamma_2}^n}\longrightarrow
D_{\gamma_2}^m \, D_{\gamma_1}^m \, w
\rightarrow \cdots \rightarrow
D_{\gamma_k}^m \, D_{\beta}^m \, w
\stackrel{D_{\gamma_k}^{-m}}\longrightarrow
D_\beta^m\, w\,.
\tag{$\star$}
\end{equation}

Every word visited by this path has the form
$D_{\gamma_i}^m 
D_{\gamma_j}^{m'} w$
for some $0\le m'\le m$
and by \lemref{lem:high-twist} has word length at least 
$m/\cc - 3d\RR \geq \RR$. Therefore, the subpath $(\star)$ stays 
outside $B_\RR$.
Also, an upper bound for $k=k(\alpha, \beta)$ 
depends on the choice 
of $\AA$ only.  Therefore the length of each subpath $(\star)$ (which is
$2km$) is of order $\RR$.
\qed

To illustrate how the switch moves work we consider the following
example.  Suppose that $\alpha-\gamma_1-\gamma_2-\beta$ is a
path in the curve complex from $\alpha$ to $\beta$ and
denote the associated Dehn twists by $a,g,h$ and $b$, respectively.  
Note that $\langle a,g \rangle$, $\langle g,h \rangle$ and $\langle h,b 
\rangle$
are free abelian subgroups of $\MCG(S)$. Let $w$ be the word
$$w=h^{-m} \ b^{m} \ g^{-m} \ h^m \ a^{-m} \ g^m.$$
After cancellation, $w$ is equivalent to $b^ma^{-m}$. Assuming that
$a, g, h$ and $b$ are all in the generating set of $\MCG(S)$, 
we can also consider $w$ as  path along the edges of the Cayley graph
of $\MCG(S)$ connecting $a^m$ to $b^{m}$ (first follow edges marked by 
the generator
$g$ for $m$ steps, then follow edges marked $a^{-1}$ for $m$ steps, and 
so on). 
Note that, in the process of carrying out the word, we stay far from the 
identity 
in the Cayley graph:  for any decomposition $w=vw'$,   
the word $w' a^m$ contains a high power of at least one
of $a, g, h$ or $b$, so its word-length is large. 
This path starts from $a^m$, travels through the quasi-flat $\langle a,g 
\rangle$ 
to $g^m$, through $\langle g,h \rangle$ to $h^m$ and through $\langle 
h,b\rangle$ 
to $b^m$.  (On the other hand, the shorter path connecting $a^m$ to 
$b^m$
corresponding to the word $b^ma^{-m}$ would go through the identity.)

\begin{figure}[ht]
\begin{center}
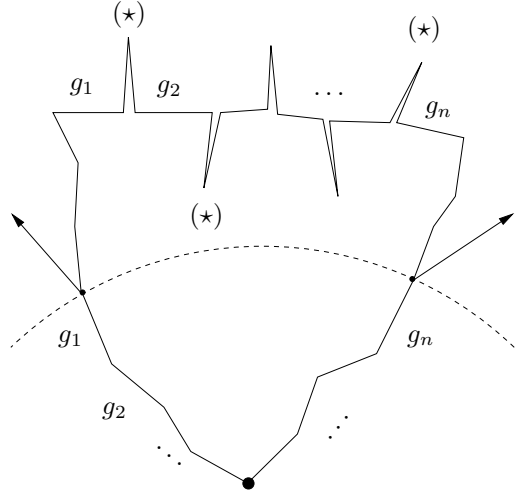
\begin{quote}
\caption{\small The path in the Cayley graph
pushes off from the ball $B_\RR$ and moves through a
chain of Dehn twist flats in $\MCG(S)$.}
\end{quote}
\end{center}
\label{fig:chainedflats}
\end{figure}

Now we will build a path $P$ from $u$ to $v$ using these
switch moves.
Let $n=|u\,v^{-1}|$ and write
$u\,v^{-1}=g_n\cdots g_2 g_1$, where
$g_i=D_{\alpha_i}^{\pm 1}$, $\alpha_i\in\AA$, is a generator.
For $0\leq r \leq n$, let  $v_r=g_r\cdots g_1\, v$,
so that
$v_0=v$ and $v_n=u$.
\begin{equation*}
\begin{split}
v
\stackrel{g_1^{\pm m}}\longrightarrow
g_1^{\pm m} v
\stackrel{g_1}\longrightarrow
g_1^{\pm m} v_1
\stackrel{(\star)}\longrightarrow
g_2^m v_1
\stackrel{g_2}\longrightarrow
g_2^mv_2
\stackrel{(\star)}\longrightarrow
g_3^mv_2
\longrightarrow \cdots \qquad \qquad \\
\cdots \longrightarrow
g_{n-1}^mv_{n-1}
\stackrel{(\star)}\longrightarrow
g_n^{\pm m} v_{n-1}
\stackrel{g_n}\longrightarrow
g_n^{\pm m} v_n
\stackrel{g_n^{\mp m}}\longrightarrow
v_n=u\,.
\end{split}
\end{equation*}

The signs $\pm$ are chosen by \lemref{lem:+or-}, which implies that the 
segments $[v,g^m_1 v]$ and
$[v,g^{-m}_1 v]$ do not both intersect $B_\RR$. We choose
the segment for the beginning of $P$ that is disjoint from $B_\RR$.
Similarly, the sign for the power at the last step is chosen so that the
path from $v_n$ to $g_n^{\pm m} v_n$ is outside $B_\RR$.
Note that $n \leq 2d\RR$ and $|v_i| \leq |v| + n \leq 3d\RR$; therefore, 
by
\lemref{lem:twist} and the assumption on $m$,
$$
d_{\alpha_j} (g_j^m v_i) \gadd m-3d\cc\RR \geq  \cc\RR\,.
$$
By \thmref{thm:wordlength}, we get  $|g_j^m v_i| \geq \RR$, which 
confirms that all of $P$ stays outside $B_\RR$.

We have shown  that the length of path $P$,
which contains $n-1$ ``switch move" subpaths, is of order $\RR^2$.
This finishes the proof of Theorem A for $\MCG(S)$.

\subsection*{Example realizing the quadratic rate.}

We now prove \thmref{thm:lower} for the case of $\MCG(S)$ by first
establishing a quasi-projection theorem.
Note 
that closest-point projection is not in general well-behaved
in the mapping class group itself, so the theorem is stated in terms of 
projection in the curve complex.
We map from $\Mod(S)$ to $\CC(S)$ by the coarse map  $g \mapsto g\AA$.  
Note that the word-length formula (\thmref{thm:wordlength}) ensures that 
this map is coarsely Lipschitz.

Distance contraction results for the mapping class group with arguments 
based on subsurface projection  
appear with  various formulations in the literature
(including Masur--Minsky, 
Behrstock, and forthcoming work of Brock--Masur--Minsky).
A statement is given and proved here in the generality which we will 
require, providing a strong parallel with \thmref{thm:QP}.

A geodesic line, ray, or segment in $\Mod(S)$ is called
{\em 
$E$-cobounded} if for every proper subsurface $Y\subset S$ and for 
every two elements $a,b\in G$, we have
$d_Y(a\AA,b\AA)\break\le E$.  

\begin{theorem}[\ (Quasi-projection for $\Mod(S)$)]
\label{thm:QP2}
For every $E>0$ there exist constants $B_1$ and $B_2$ depending on $E$ 
and the topology of $S$ such that the following holds.
Let $G\subset \Mod(S)$ be an 
$E$-cobounded geodesic, suppose $g\in\Mod(S)$ satisfies $d(g,G)>B_1$,
and let $R={d(g,G)}/{B_1}$.
Then $\G=G\AA$ is a quasi-geodesic in $\CC(S)$ and  
$$
\diam\bigl(\Proj_\G( B_R(g)\AA)\big)\le B_2\,.
$$
\end{theorem}

As before, the theorem could be stated for 
quasi-geodesics with the same 
argument.
Note that for quasi-projection to a segment, it suffices that the 
endpoints $a,b\in \Mod(S)$ satisfy $d_Y(a\AA,b\AA)\le E$ 
in order for the conclusion to obtain.  

For a concrete application, consider a pseudo-Anosov element $w$ and its 
axis $\{w^n\}$ in $\Mod(S)$.
These axes are known to be cobounded by work of
Masur--Minsky.  We then 
apply the quasi-projection theorem
exactly as above to show that, for a high power $m$, 
any path $P$ connecting $w^{-m}$ to $w^m$ outside the ball of radius 
$t:=|w^m|\asymp |m|$
must have length at least on the order of $t^2$.  
(Since projection to 
the curve complex coarsely contracts distances, this
agrees with the desired inequality.)

The rest of this section is devoted to proving \thmref{thm:QP2}.
The constants in the rest of this section will depend on the choice of 
$E$,
which is chosen once and for all. Keeping this in mind we can write,
for example, 
$$
d_S(\AA, \G) \asymp 1\,.
$$

\PR{Proof of \thmref{thm:QP2}}
Let $h\in B_R(g)$.
Let $\alpha=g(\AA)$ and $\beta = h(\AA)$ (these 
are sets of curves, and they fill $S$).
Let $\ab$ and $\bb$ be closest-point projections of
$\alpha$ and $\beta$ \hbox{to $\G$}, respectively.
Choose $\gb, \hb \in G$ so that
$$
\gb (\AA)\cap \ab \neq\emptyset  \quad \text{and}\quad \hb (\AA)\cap \bb 
\neq \emptyset\,, 
$$
which is possible because $\G=G\AA$.
Our goal is to show that $d_S(\ab,\bb)\le B_2$, so assume for 
contradiction that 
$$d_S ( \ab , \bb ) > B_2\,.
$$

\begin{figure}[ht]
\begin{center}
\input{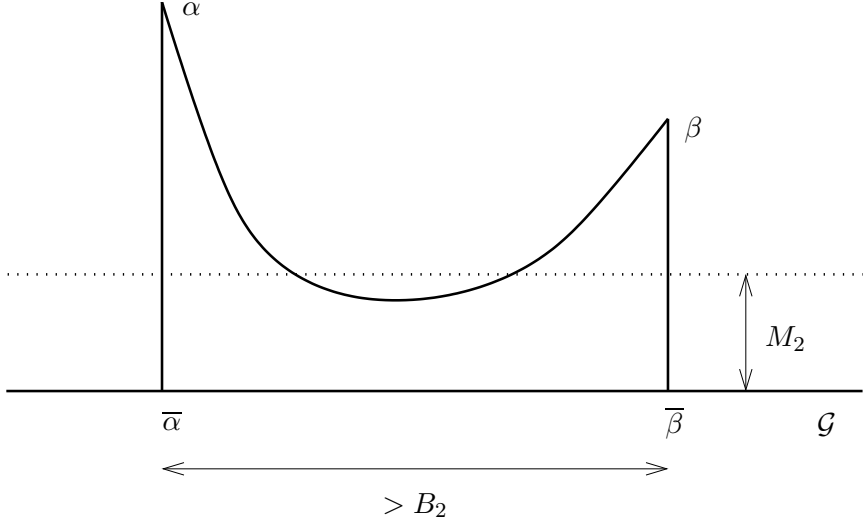}
\begin{quote}
\caption{\small We assume for contradiction that the
projections of $\alpha$ 
and 
$\beta$ to $\G$ are far apart in $\CC(S)$.}
\end{quote}
\end{center}
\label{Fig:projection}
\end{figure}

Given any $\thresh,\cc$ as in \thmref{thm:wordlength}, we can define
 $\cY$ to be the set of proper subsurfaces $Y \subsetneq S$ with 
$d_Y(\alpha, \ab) \geq \thresh$.  By the theorem's hypothesis, we have
$|g^{-1}\gb|=d(g,\gb)\ge d(g,G)= B_1R$ and so 
$$
\sum_{Y\in \cY} d_Y(\alpha, \ab)   + d_S(\alpha, \ab) 
\geq \frac {|g^{-1} \gb|}{ \cc} \geq \frac{B_1R}{\cc}\,.
$$
We will proceed in two cases. First assume that
\begin{equation} \label{eq:condition}
d_S(\alpha, \ab) \geq \frac{B_1 R}{2\cc}\,.
\end{equation}
If $B_2 > M_1$, then any geodesic connecting $\alpha$ to
$\beta$ passes through the $M_2$-neigh\-bor\-hood of
$\G$, by 
\lemref{lem:quad}.
Hence for large enough $B_1$, we have
$$
d_S(\alpha, \beta) \geq d_S(\alpha, \ab) - M_2 \geq 
 \frac{B_1 R}{2\cc} - M_2 \geq  \frac{B_1 R}{3\cc}\,.
$$
Enlarging $B_1$ again if necessary, we can assume $d_S(\alpha, \beta) 
\geq \thresh_0$.
\thmref{thm:wordlength} implies
$$
|h^{-1}g|\geq\frac 1{\cc}d_S(\alpha,\beta)\geq\frac{B_1R}
{3\cc 
\cc_0}
> R\,,
$$
for  $B_1$ sufficiently large.  But this contradicts the assumption that 
$d(g,h)\le R$.

In the second case, the assumption of Equation \eqref{eq:condition} is 
not true, 
so we have 
\begin{equation}\label{eq:cond'}
\sum_{Y\in \cY} d_Y(\alpha, \ab) \geq \frac{B_1 R}{2\cc},
\end{equation}
where again $\cY$ is chosen relative to the threshold $\thresh$.

Next we show that the subsurfaces with large projection
distance between $\alpha$ and $\ab$ do not also have a
large distance between $\beta$ and 
$\bb$.

\begin{thm}{Claim}
For every $Y \in \cY$, $d_Y(\beta, \bb) \leq M_0$.
\end{thm}

\pr 
If $\thresh > M_0$, then \thmref{thm:sub-proj} implies that
any geodesic $[\alpha,\ab]$ intersects the
$1$-neighborhood of $\p Y$
for every $Y\in\cY$.
Applying \lemref{lem:shadow} and \lemref{lem:Lip} 
we can conclude that
$$
d_S\big(\ab, \Proj(\partial Y)\big) \leq M_3 + 2 M_4\,. 
$$
If $d_Y(\beta, \bb) > M_0$ for any particular $Y\in \cY$, we similarly 
get
$$
d_S(\bb, \Proj\big(\partial Y)\big) \leq M_3 + 2 M_4\,. 
$$
But then
$$
d_S(\ab, \bb) \leq 2M_3 + 4 M_4\,.
$$
For $B_2$ larger than this final constant, this is a contradiction.
\qed

Note that $\gb,\hb\in G$ and $G$ is a $E$-cobounded
geodesic. 
Therefore, for every $Y \in \cY$,
$$
d_Y(\ab,\bb)=d_Y\big(\gb(\AA),\hb(\AA)\big)\leq E\,.
$$
Using the triangle inequality for the projection distance 
(\lemref{lem:triangle-inequality}) we get
$$ 
d_Y(\alpha, \beta) \gadd
d_Y(\alpha, \ab) -d_Y(\ab, \bb)- d_Y(\bb, \beta)
\geq d_Y(\alpha, \ab) - M_0 - E\,.
$$

By choosing $\thresh$ large enough, we can ensure for
$Y\in\cY$
that $d_Y(\alpha, \beta)$ is
much larger than these additive errors in order to get a 
multiplicatively 
coarse equality
$$
d_Y(\alpha, \beta) \emul d_Y(\alpha, \ab)
$$
as well as ensuring that
$$
d_Y(\alpha, \beta) \geq \thresh_0\,.
$$
We apply these, as well as \thmref{thm:wordlength} one last time, to see 
that 
\begin{align*}
\cc_0 \, |h^{-1}g|&\gadd\sum_{Y\in\cY}\big[d_Y(g(\AA),h(
\AA))\big]_{\thresh_0} 
\\
&\eadd\sum_{Y\in\cY}\big[d_Y(\alpha,\beta)\big]_{\thresh_
0} \\
&\emul\sum_{Y\in\cY}d_Y(\alpha,\ab)\geq\frac{B_1R}{2\cc}
\,, 
\end{align*}
where the final inequality comes from \eqref{eq:cond'}.
We have shown that 
$$
 |h^{-1}g|  \succ \frac{B_1 R}{2\cc\cc_0}\,.
$$
Again, choosing $B_1$ large enough provides the
contradiction.
\qed

\bigskip\bigskip

\noindent{\sc Moon Duchin}, 
Department of Mathematics, University of Michigan, Ann
Arbor, MI 48109, USA
\hfill{\tt mduchin@umich.edu}

\rightline{\tt http://www.math.lsa.umich.edu/$\sim$mduchin/}

\medskip
\noindent{\sc Kasra Rafi}, 
Department of Mathematics, University of Oklahoma,
Norman, OK 73019, USA \hfill{\tt rafi@math.ou.edu}

\rightline{\tt http://www.math.ou.edu/$\sim$rafi/}

\bigskip

\rightline{Received: April 5, 2008}
\rightline{Revision: October 16, 2008}
\rightline{Accepted: November 23, 2008}


\small

\begin{thebibliography}{XXX}
\itemsep=0pt
\parskip=1pt

\bibitem[A]{abikoff:TT}
{\sc W.\ Abikoff},
{The Real Analytic Theory of {T}eichm\"uller Space},
Springer Lecture Notes in Mathematics 820,
  Springer, Berlin (1980).

\bibitem[B]{behrstock:AG}
{\sc J.\ Behrstock},
  Asymptotic geometry of the mapping class group and 
{T}eichm\"uller
  space,
  {Geometry \& Topology} 10 (2006), 1523--1578.

\bibitem[Be]{bers:EP}
{\sc L.\ Bers},
  An extremal problem for quasiconformal mappings and a theorem 
by
  {T}hurston,
  {Acta Math.} 141:1-2 (1978), 73--98.

\bibitem[BrH]{bridson:NPC}
{\sc M.R.\ Bridson, A.\ Haefliger},
{Metric Spaces of Non-positive Curvature}, 
Springer Grundlehren der Mathematischen Wissenschaften
[Fundamental Principles 
of
  Mathematical Sciences] 319 (1999).

\bibitem[G1]{gersten:QD}
{\sc S.M.\ Gersten},
Quadratic divergence of geodesics in {${\rm CAT}(0)$}
spaces,
{Geom.\ Funct.\ Anal.} 4:1 (1994), 37--51.

\bibitem[G2]{gersten:D3}
{\sc S.M.\ Gersten},
  Divergence in $3$-manifold groups,  
  {Geom.\ Funct.\ Anal.}  4:6 (1990), 633--647.

\bibitem[Gr]{gromov:AI}
{\sc M.\ Gromov},
{Geometric Group Theory, Vol.\ 2: Asymptotic Invariants of
Infinite Groups}
(Niblo, Roller, eds), LMS Lecture Note Ser. 182,
Cambridge Univ.\ Press, Cambridge (1993).

\bibitem[H]{hubbard:TT}
{\sc J.\ Hubbard},
  {Teichm\"uller Theory and Applications to Geometry, 
Topology and
  Dynamics},
  Matrix Editions, Ithaca, NY, 2006.

\bibitem[I]{ivanov:ACC}
{\sc N.V.\ Ivanov},
  Automorphism of complexes of curves and of {T}eichm\"uller 
spaces, {Internat.\ Math.\ Res.\ Notices} 1997:14 (1997),
651--666.

\bibitem[KL]{kapovich:3M}
{\sc M.\ Kapovich, B.\ Leeb},
  $3$-manifold groups and nonpositive curvature,
  {Geom.\ Funct.\ Anal.}  8:5 (1998), 841--852.

\bibitem[Ke]{kerckhoff:AG}
{\sc S.P.\ Kerckhoff},
  The asymptotic geometry of {T}eichm\"uller space,
  {Topology} 19:1(1980), 23--41.

\bibitem[M1]{masur:CG}
{\sc H.A.\ Masur},
  On a class of geodesics in {T}eichm\"uller space,
  {Ann.\ of Math. (2)} 102:2 (1975), 205--221.

\bibitem[M2]{masur:TP}
{\sc H.A.\ Masur},
  Transitivity properties of the horocyclic and geodesic flows 
on moduli space,
  {J.\ Analyse Math.} 39 (1981), 1--10.

\bibitem[MM1]{minsky:CCI}
{\sc H.A.\ Masur, Y.N.\ Minsky},
  Geometry of the complex of curves. {I}. {H}yperbolicity,
{Invent.\ Math.} 138:1 (1999), 103--149.

\bibitem[MM2]{minsky:CCII}
{\sc H.A.\ Masur, Y.N.\ Minsky},
  Geometry of the complex of curves. {II}. {H}ierarchical 
structure,
  {Geom.\ Funct.\ Anal.} 10:4 (2000), 902--974.

\bibitem[MW]{masur:NH}
{\sc H.A.\ Masur, M.\ Wolf},
  Teichm\"uller space is not {G}romov hyperbolic,
{Ann.\ Acad.\ Sci.\ Fenn.\ Ser.\ A I Math.} 20:2 (1995), 
259--267. 

\bibitem[Mi1]{minsky:PR}
{\sc Y.N.\ Minsky},
  Extremal length estimates and product regions in 
{T}eichm\"uller
  space,
  {Duke Math.\ J.} 83:2(1996), 249--286.

\bibitem[Mi2]{minsky:QP}
{\sc Y.N.\ Minsky},
  Quasi-projections in {T}eichm\"uller space,
  {J.\ Reine Angew.\ Math.} 473 (1996), 121--136.

\bibitem[P]{papasoglu:SQ}
{\sc P.\ Papasoglu},
  On the subquadratic isoperimetric inequality,
in ``{Geometric Group Theory}",
de Gruyter (1995), 149--157.

\bibitem[R1]{rafi:CM}
{\sc K.\ Rafi},
  {A combinatorial model for the Teichm\"uller metric},
 {Geom.\ Funct.\ Anal.} 17:3 (2007), 936--959.

\bibitem[R2]{rafi:thesis}
{\sc K.\ Rafi},
 {A characterization of short curves of a geodesic in 
Teichm\"uller space},
  {Geometry \& Topology} 9 (2005), 179--202.

\bibitem[S]{short:GGT}
{\sc H.\ Short, ed.},
  MSRI notes on hyperbolic groups,
in ``Group Theory from a Geometrical Viewpoint",  
World Sci.\ Publ., River Edge, NJ (1991), 3--63. 

\bibitem[St]{strebel:QD}
{\sc K.\ Strebel},
{Quadratic Differentials}, vol.\,5, of {A Series of 
Modern
  Surveys in Mathematics},
  Springer-Verlag, Berlin, 1980.

\bibitem[W1]{wolpert:NC}
{\sc S.A.\ Wolpert},
  Noncompleteness of the {W}eil-{P}etersson metric for 
{T}eichm\"uller
  space,
  {Pacific J.\ Math.} 61:2 (1975), 573--577.

\bibitem[W2]{wolpert:LS}
{\sc S.A.\ Wolpert},
The length spectra as moduli for compact Riemann surfaces,
{Ann.\ of Math. (2)} 109:2 (1979), 323--351.

\end{thebibliography}
\end{document}